\documentclass[preprint,12pt]{elsarticle}
\usepackage{amsmath,amssymb,amsfonts}
\usepackage{amsthm}

\usepackage[colorlinks=true,linkcolor=blue!50!black,citecolor=blue!50!black,urlcolor=blue!50!black]{hyperref}

\newtheorem{theorem}{Theorem}

\newtheorem{lemma}[theorem]{Lemma}
\newtheorem{proposition}[theorem]{Proposition}

\newcommand{\Z}{\mathbb{Z}}

\journal{European Journal of Combinatorics}

\begin{document}

\begin{frontmatter}



\title{Positivity of Signed Odd-Multiplicity Partitions}


\author{Mircea Merca}

\affiliation{organization={Department of Mathematical Methods and Models, Fundamental Sciences Applied in Engineering Research Center, National University of Science and Technology Politehnica Bucharest},
            city={Bucharest},
            postcode={060042}, 
            country={Romania}}

\affiliation{organization={Academy of Romanian Scientists},
	city={Bucharest},
	postcode={050044}, 
	country={Romania}}

\begin{abstract}
	We establish coefficientwise nonnegativity for the signed
	partition $q$-series
	\[
	A_k(q)=
	\sum_{1\le j_1<\cdots<j_k}
	\prod_{i=1}^k\frac{q^{j_i}}{1+q^{2j_i}}
	=\sum_{n\ge0}a_k(n)q^n,
	\qquad k\ge0,
	\]
	with $A_0(q)=1$.
	The coefficient $a_k(n)$ counts, with signs, partitions
	of $n$ having exactly $k$ distinct part sizes, each
	occurring with odd multiplicity. A partition with
	multiplicities $2t_1+1,\ldots,2t_k+1$ is assigned the
	sign $(-1)^{t_1+\cdots+t_k}$.
	We prove that $a_k(n)\ge0$ for all $k,n\ge0$,
	thereby resolving a conjecture on signed
	odd-multiplicity partitions.
	Our proof combines a positive theta--MacMahon
	factorization for the generating functions restricted
	to odd part sizes with a dyadic product decomposition
	that extends positivity to arbitrary part sizes.
	Moreover, we establish strict positivity
	$a_k(n)>0$ for every $k\ge2$ and
	$n\ge\binom{k+1}{2}$.
\end{abstract}




\end{frontmatter}



\section{Introduction}\label{sec:introduction}

For each nonnegative integer $k$, consider the following signed
partition $q$-series family
\begin{equation}\label{eq:Adef}
A_k(q)=
\sum_{1\le\lambda_1<\cdots<\lambda_k}
\prod_{i=1}^k
\frac{q^{\lambda_i}}{1+q^{2\lambda_i}}
=\sum_{n\ge0}a_k(n)\,q^n,
\end{equation}
with $A_0(q)=1$. The coefficients of the series $A_k(q)$ have a natural interpretation as signed counts of partitions.

\smallskip

For integers $k,n\ge0$, let
$\mathcal P_k^{\mathrm{odd}}(n)$ denote the set of
partitions of $n$ having exactly $k$ distinct part
sizes, each occurring with odd multiplicity.
We adopt the convention
$\mathcal P_0^{\mathrm{odd}}(0)=\{\varnothing\}$
and $\mathcal P_0^{\mathrm{odd}}(n)=\varnothing$
for $n>0$.
For $k\ge1$, every partition $\pi$ in
$\mathcal P_k^{\mathrm{odd}}(n)$ can be written uniquely as
\begin{equation}\label{eq:partition-form}
n=(2t_1+1)\lambda_1+\cdots+
(2t_k+1)\lambda_k,
\end{equation}
where $1\le\lambda_1<\cdots<\lambda_k$ and
$t_1,\ldots,t_k$ are nonnegative integers.
Assign to $\pi$ the sign
$(-1)^{t_1+\cdots+t_k}$, with the empty
partition assigned sign $1$.
Since
\[
\frac{q^j}{1+q^{2j}}
=\sum_{t\ge0}(-1)^tq^{(2t+1)j},
\]
the coefficients in \eqref{eq:Adef} satisfy
\begin{equation}\label{eq:signedcount}
a_k(n)=
\sum_{\pi\in\mathcal P_k^{\mathrm{odd}}(n)}
(-1)^{t_1+\cdots+t_k}.
\end{equation}

Let $N_k^+(n)$ and $N_k^-(n)$ denote the numbers of
partitions in $\mathcal P_k^{\mathrm{odd}}(n)$ with
positive and negative signs, respectively. Then
\begin{equation}\label{eq:signed-difference}
a_k(n)=N_k^+(n)-N_k^-(n).
\end{equation}
Cancellations occur even for small values of $k$.
For example, when $k=3$ and $n=8$, the partitions
$5+2+1$ and $4+3+1$ have positive sign, whereas
$3+2+1+1+1$ has negative sign. Thus $a_3(8)=1$.
The sign of \eqref{eq:signed-difference} is therefore
not apparent from its combinatorial definition.

\smallskip
The series $A_k(q)$ coincide with the
specialization $U_k(0,q)$ studied by
Amdeberhan, Andrews, and Tauraso
\cite{AATFurther}. Their interpretation
as signed partition enumerators was
developed in \cite{MercaSigned}, where
coefficientwise nonnegativity was conjectured.
The case $k=2$ was subsequently proved
in \cite{MercaTwo}.

\smallskip

A further connection arises from the partition
function $\operatorname{pod}(n)$, which counts
partitions of $n$ into distinct odd parts and
unrestricted even parts; see
\cite{HirschhornSellersPod}.
Its generating function is
\begin{equation}\label{eq:pod}
\sum_{n\ge0}\operatorname{pod}(n)q^n
=\frac{(q^2;q^4)_\infty}{(q;q)_\infty}
=\frac{1}{\psi(-q)},
\end{equation}
where
\[
\psi(q)=\sum_{r\ge0}q^{r(r+1)/2},
\qquad
(a;q)_\infty=\prod_{j\ge0}(1-aq^j).
\]
%
%
In \cite[Theorem~1.5]{MercaSigned}, it was shown
that
\begin{equation}\label{eq:partial-theta}
A_k(q)=
\frac{q^{\binom{k+1}{2}}}{\psi(-q)}
\Theta_k^*(q),
\end{equation}
where
\begin{equation}\label{eq:theta-star}
\Theta_k^*(q)=
\sum_{r\in\mathbb Z}
(-1)^r\binom{r+k}{k}
q^{r(2r+2k+1)}.
\end{equation}
Here the binomial coefficients are interpreted
in the polynomial sense:
\[
\binom{r+k}{k}
=\frac{(r+1)(r+2)\cdots(r+k)}{k!},
\qquad r\in\mathbb Z.
\]
The equivalence with the unilateral representation
in \cite{MercaSigned} follows by pairing the terms
indexed by $r$ and $-r-k-1$.

Combining \eqref{eq:partial-theta} with
\eqref{eq:pod}, we obtain
\begin{equation*}
a_k\left(n+\binom{k+1}{2}\right)
=
\sum_{r\in\mathbb Z}
(-1)^r\binom{r+k}{k}
\operatorname{pod}
\bigl(n-r(2r+2k+1)\bigr),
\end{equation*}
for all $k,n\ge0$, with the convention
$\operatorname{pod}(m)=0$ for $m<0$.
Thus the conjectured nonnegativity of
$a_k(n)$ is equivalent to the nonnegativity
of these alternating weighted sums, a property
that is far from evident from their definition.

\smallskip
Our main result is the following.

\begin{theorem}[Positivity of signed
	odd-multiplicity partitions]\label{thm:positive}
	For all integers $k,n\ge0$, the coefficients of
	$A_k(q)$ satisfy $a_k(n)\ge0$.
	Moreover, for every $k\ge2$, the coefficients of $A_k(q)$ vanish below degree $\binom{k+1}{2}$ and satisfy $a_k(n)>0$ for $n\ge\binom{k+1}{2}$. Equivalently, for every $k\ge2$, the series
	\begin{equation}\label{eq:positive-partial-product}
	\frac{(q^2;q^4)_\infty}{(q;q)_\infty} \sum_{r\in\mathbb{Z}}(-1)^r\binom{r+k}{k} q^{r(2r+2k+1)} 
	\end{equation}
has strictly positive integer coefficients.
\end{theorem}

Our proof proceeds in two steps. We first consider
the generating functions obtained by restricting
all part sizes to odd integers:
\begin{equation}\label{eq:Bdef}
B_j(q)=
\sum_{\substack{1\le s_1<\cdots<s_j\\
		s_i\ \mathrm{odd}}}
\prod_{i=1}^j
\frac{q^{s_i}}{1+q^{2s_i}},
\qquad B_0(q)=1.
\end{equation}
Sellers and Tauraso \cite{ST} established transfer
identities relating these series to MacMahon-type
sums and theta--Lambert series. While these
identities provide the starting point for our
argument, positivity is not immediately apparent
in the odd-indexed case. We obtain a positive
factorization using the Jacobi triple product,
thereby proving coefficientwise nonnegativity
for every $B_j(q)$.

\smallskip
The second step extends this positivity to
arbitrary part sizes. Using the unique
decomposition of each positive integer as
$2^\nu s$, where $\nu\ge0$ and $s$ is odd,
we establish the dyadic factorization
\begin{equation}\label{eq:dyadic-intro}
\sum_{k\ge0}A_k(q)z^k
=
\prod_{\nu\ge0}
\left(
\sum_{j\ge0}B_j(q^{2^\nu})z^j
\right).
\end{equation}
The combination of the positive odd-part-size
factorization and this dyadic decomposition
proves the conjectured nonnegativity for all
$k$. A further analysis of the positive
factors yields strict positivity throughout
the admissible range for $k\ge2$.

\smallskip
All identities are interpreted as formal
power-series identities. We write
$F(q)\succeq0$ when
$F(q)\in\mathbb Z[[q]]$ has nonnegative
coefficients, and $[q^n]F$ denotes the
coefficient of $q^n$ in $F(q)$. Section~\ref{sec:transfer}
establishes the positive odd-size identities,
Section~\ref{sec:dyadic} lifts them to arbitrary part sizes,
and Section~\ref{sec:strict} proves strict positivity.

\section{The odd-part-size series and the known MacMahon transfer}
\label{sec:transfer}

Recall the odd-part-size series $B_j(q)$ defined in
\eqref{eq:Bdef}.  In the notation of \cite{ST},
$B_j(q)=\widetilde U_j(0,q)$.
For $m\ge0$, put
\begin{align}
M_m(q)&=\sum_{1\le r_1<\cdots<r_m}
\prod_{i=1}^{m}\frac{q^{r_i}}{(1-q^{r_i})^2},\label{eq:Mdef}\\
H_m(q)&=\sum_{\substack{1\le s_1<\cdots<s_m\\s_i\ \mathrm{odd}}}
\prod_{i=1}^{m}\frac{q^{s_i}}{(1-q^{s_i})^2},\label{eq:Hdef}
\end{align}
where $M_0=H_0=1$.  These are MacMahon's ordinary and odd-part
sums of divisors.  They have nonnegative integer coefficients because
\begin{equation}\label{eq:positive-factor}
\frac{q^r}{(1-q^r)^2}=\sum_{h\ge1}h\,q^{hr}.
\end{equation}

\begin{proposition}[Odd-size transfers and positive Jacobi factorization]\label{prop:transfer}
	For every $m\ge0$,
	\begin{equation}\label{eq:transfers}
	B_{2m}(q)=H_m(q^4),\qquad
	B_{2m+1}(q)=q\,\psi(q^4)^2\,M_m(q^8).
	\end{equation}
	Consequently $B_j(q)\succeq0$ for every $j\ge0$.
\end{proposition}

\begin{proof}
	Consider the generating function
	\begin{equation}\label{eq:G}
	G(z,q)=\sum_{j\ge0}B_j(q)z^j
	=\prod_{r\ge1}
	\frac{1+zq^{2r-1}+q^{4r-2}}
	{1+q^{4r-2}}.
	\end{equation}
	Set $z=u+u^{-1}$ and introduce
	\[
	\Theta(x)=\sum_{h\in\mathbb Z}(-1)^h x^{h^2}
	=(x^2;x^2)_\infty(x;x^2)_\infty^2.
	\]
	By Jacobi's triple product identity
\cite[Theorem~2.8]{AndrewsPartitions},
	\[
	\sum_{h\in\mathbb Z}u^h q^{h^2}
	=(q^2;q^2)_\infty
	\prod_{r\ge1}
	(1+uq^{2r-1})(1+u^{-1}q^{2r-1}).
	\]
	Since
	\[
	(q^2;q^2)_\infty
	\prod_{r\ge1}(1+q^{4r-2})
	=\Theta(q^4),
	\]
	equation~\eqref{eq:G} becomes
	\begin{equation}\label{eq:thetaquotient}
	G(z,q)=
	\frac{1}{\Theta(q^4)}
	\sum_{h\in\mathbb Z}u^h q^{h^2}.
	\end{equation}
	
	Put $x=q^4$. Separating the even and odd indices
	in the numerator of \eqref{eq:thetaquotient},
	the ordinary and shifted forms of the Jacobi
	triple product give
	\begin{align*}
	\sum_{h\in\mathbb Z}u^{2h}q^{4h^2}
	&=(x^2;x^2)_\infty
	\prod_{r\ge1}
	(1+u^2x^{2r-1})(1+u^{-2}x^{2r-1}),\\
	\sum_{h\in\mathbb Z}u^{2h+1}q^{(2h+1)^2}
	&=qz(x^2;x^2)_\infty
	\prod_{r\ge1}
	(1+u^2x^{2r})(1+u^{-2}x^{2r}).
	\end{align*}
	Here the factor $z$ in the second identity
	comes from
	$u(1+u^{-2})=u+u^{-1}=z$.
	
	Using
	\[
	(1+u^2t)(1+u^{-2}t)
	=(1-t)^2+z^2t
	\]
	together with
	\[
	\Theta(x)
	=(x^2;x^2)_\infty
	\prod_{r\ge1}(1-x^{2r-1})^2
	\]
	and
	\[
	\frac{(x^2;x^2)_\infty^3}{\Theta(x)}
	=\psi(x)^2,
	\]
	we obtain the following expressions for the
	even and odd parts of $G(z,q)$:
	\begin{align}
	G_{\rm e}(z,q)
	&=\prod_{\substack{s\ge1\\s\ {\rm odd}}}
	\left(
	1+\frac{z^2q^{4s}}{(1-q^{4s})^2}
	\right),
	\label{eq:positive-even}\\
	G_{\rm o}(z,q)
	&=zq\psi(q^4)^2
	\prod_{r\ge1}
	\left(
	1+\frac{z^2q^{8r}}{(1-q^{8r})^2}
	\right).
	\label{eq:positive-odd}
	\end{align}
	These identities hold formally: for each fixed
	power of $q$, the corresponding symmetric Laurent
	polynomial in $u$ is a polynomial in
	$z=u+u^{-1}$, and its parity is preserved.
	
	Extracting the coefficients of $z^{2m}$
	and $z^{2m+1}$ yields, respectively,
	\[
	B_{2m}(q)=H_m(q^4),
	\qquad
	B_{2m+1}(q)=q\psi(q^4)^2M_m(q^8),
	\]
	which proves \eqref{eq:transfers}.
	
	The even-index identity agrees with
	\cite[Eq.~(12)]{ST}. For odd indices,
	Sellers and Tauraso \cite[Eqs.~(13)--(14)]{ST}
	obtained the transfer
	\[
	B_{2m+1}(q)=qW_m(q^4),
	\]
	where $W_m$ is defined by their
	theta--Lambert series.
	Comparison with \eqref{eq:transfers} therefore
	gives the positive factorization
	\begin{equation}\label{eq:Wpositive}
	W_m(q)=\psi(q)^2M_m(q^2).
	\end{equation}
	
	Finally, $H_m(q)$ and $M_m(q)$ have
	nonnegative integer coefficients by
	\eqref{eq:positive-factor}, as does $\psi(q)$.
	Consequently, $B_j(q)\succeq0$ for every
	$j\ge0$.
\end{proof}

\section{Dyadic factorization and coefficientwise positivity}
\label{sec:dyadic}
Every positive integer has a unique decomposition $2^\nu s$
with $\nu\ge0$ and $s$ odd.  This elementary observation is the
key to transferring positivity from $B_j$ to every $A_k$.

\begin{proposition}[Positive dyadic factorization]\label{prop:dyadic}
	With $G(z,q)$ defined by \eqref{eq:G}, one has
	\begin{equation}\label{eq:dyadic-factorization}
	\sum_{k\ge0}A_k(q)\,z^k
	=\prod_{\nu\ge0}G(z,q^{2^\nu}).
	\end{equation}
	In particular, with $G_{\rm e},G_{\rm o}$ defined in
	\eqref{eq:positive-even}--\eqref{eq:positive-odd},
	\begin{equation}\label{eq:positive-full}
	\sum_{k\ge0}A_k(q)\,z^k
	=\prod_{\nu\ge0}
	\bigl(G_{\rm e}(z,q^{2^\nu})+G_{\rm o}(z,q^{2^\nu})\bigr).
	\end{equation}
	All coefficients of the right-hand side are nonnegative integers.
	Consequently, $a_k(n)\ge0$ for every $k,n\ge0$.
\end{proposition}
\begin{proof}
	Introduce the elementary-symmetric product
	\begin{equation}\label{eq:all-prod}
	\sum_{k\ge0}A_k(q)\,z^k
	=\prod_{r\ge1}\left(1+z\frac{q^r}{1+q^{2r}}\right).
	\end{equation}
	For each $r=2^\nu s$ with $s$ odd, its factor is
	\[
	1+z\frac{(q^{2^\nu})^s}{1+(q^{2^\nu})^{2s}}.
	\]
	Regrouping the factors of \eqref{eq:all-prod} according to $\nu$
	proves \eqref{eq:dyadic-factorization}; inserting the
	positive even--odd decomposition of $G$ proves
	\eqref{eq:positive-full}.
	These products are well-defined in $\Z[z][[q]]$: the factor for
	level $\nu$ is $1+O(q^{2^\nu})$, so a fixed coefficient of
	$q^N$ depends on only finitely many levels and part sizes.
	Each factor of $G_{\rm e}$ has nonnegative coefficients by
	\eqref{eq:positive-factor}.  The same holds for $G_{\rm o}$,
	since $\psi(q^4)^2$ has nonnegative integer coefficients.
	The coefficient of $z^kq^n$ in \eqref{eq:positive-full}
	is therefore a finite sum of nonnegative integers.
	By \eqref{eq:all-prod} it equals $a_k(n)$.  This proves
	the nonnegativity assertion of Theorem~\ref{thm:positive}.
\end{proof}

Splitting the part sizes into odd and even gives the useful recurrence
\begin{equation}\label{eq:recurrence}
A_k(q)=\sum_{j=0}^k B_j(q)A_{k-j}(q^2).
\end{equation}

\section{Strict positivity}\label{sec:strict}
Let $N_k=k(k+1)/2$.  By \eqref{eq:Adef}, the least possible
degree of $A_k(q)$ is $N_k$.  Moreover, for every $k\ge1$,
\begin{equation}\label{eq:two-initial}
a_k(N_k)=a_k(N_k+1)=1.
\end{equation}
Indeed, increasing an odd multiplicity from one increases the
weight by at least two; the unique increasing $k$-tuples with
sums $N_k$ and $N_k+1$ are $(1,\ldots,k)$ and
$(1,\ldots,k-1,k+1)$, respectively.

The following observation on odd part sizes is sufficient; no
quantitative coefficient estimate is needed.
\begin{lemma}\label{lem:positiveB}
	For every integer $j\ge2$ and $r\ge0$,
	\[
	[q^{j^2+4r}]B_j(q)>0.
	\]
\end{lemma}
\begin{proof}
	If $j=2m$, select the smallest $m$ odd sizes
	$1,3,\ldots,2m-1$ in $H_m(q^4)$.  Their contribution is
	\[
	\frac{q^{j^2}}{\prod_{i=1}^m(1-q^{4(2i-1)})^2}.
	\]
	Its factor $(1-q^4)^{-2}$ has a positive coefficient at every
	power $q^{4r}$.
	
	If $j=2m+1$, where $m\ge1$, retain the smallest $m$ sizes
	$1,\ldots,m$ in $M_m(q^8)$.  Since
	$1+8\binom{m+1}{2}=j^2$, the identity
	\eqref{eq:transfers} contains the nonnegative series
	\[
	\frac{q^{j^2}\psi(q^4)^2}
	{\prod_{i=1}^m(1-q^{8i})^2}.
	\]
	The constant and $q^4$ coefficients of $\psi(q^4)^2$ are
	positive, while $(1-q^8)^{-2}$ has a positive coefficient at
	every power $q^{8r}$.  These two residue classes cover all
	exponents $j^2+4r$.
\end{proof}

\begin{proof}[Proof of Theorem~\ref{thm:positive}]
	The vanishing below $N_k$ follows from the defining series.
	The case $k=2$ was established in \cite{MercaTwo}.
	
	For $k=3$, the recurrence \eqref{eq:recurrence} contains
	$B_2(q)A_1(q^2)$ and $B_1(q)A_2(q^2)$, both nonnegative.
	Since $a_1(1)=a_1(2)=1$, Lemma~\ref{lem:positiveB} shows
	that the first product has positive coefficients at every
	even degree $n\ge6$.  In every odd degree $n\ge7$, the
	second product contributes at least
	\[
	[q]B_1(q)\,a_2\bigl((n-1)/2\bigr)
	=a_2\bigl((n-1)/2\bigr)>0.
	\]
	
	Now let $k\ge4$.  Write $k=2m$ or $k=2m+1$ and choose
	\[
	\begin{array}{c|cc}
	& (j_0,s_0)&(j_1,s_1)\\ \hline
	k=2m & (m,m)&(m+1,m-1)\\
	k=2m+1 & (m+1,m)&(m,m+1)
	\end{array}
	\]
	so that $j_i+s_i=k$, $j_i\ge2$, $s_i\ge1$, and
	\begin{equation}\label{eq:degrees}
	j_0^2+2N_{s_0}=N_k,\qquad
	j_1^2+2N_{s_1}=N_k+1.
	\end{equation}
	Both $B_{j_0}(q)A_{s_0}(q^2)$ and
	$B_{j_1}(q)A_{s_1}(q^2)$ occur in \eqref{eq:recurrence}
	with nonnegative coefficients.  Lemma~\ref{lem:positiveB}
	and \eqref{eq:two-initial} show that the first product has
	positive coefficients at degrees $N_k+4r$ and $N_k+2+4r$,
	whereas the second has positive coefficients at
	$N_k+1+4r$ and $N_k+3+4r$, for every $r\ge0$.
	Indeed, every integer $n\ge N_k$ has a unique representation
	$n=N_k+\delta+4r$ with $\delta\in\{0,1,2,3\}$ and $r\ge0$;
	these are exactly the four progressions just covered.
	Thus $a_k(n)>0$ for every $n\ge N_k$.
\end{proof}

\section{Concluding remark}
The positive theta--MacMahon factorization for odd part sizes,
combined with dyadic decomposition, proves the conjectured
coefficientwise nonnegativity of the full family
$A_k(q)=U_k(0,q)$.  The same framework establishes strict
positivity for every $k\ge2$ from the minimal degree
$\binom{k+1}{2}$ onward.

%
%
%




  \bibliographystyle{elsarticle-num-names} 


\end{document}